\documentclass[11pt,a4paper]{article}

\usepackage{amsmath}
\usepackage{amsthm}
\usepackage{amsfonts}

\parskip\smallskipamount
\newtheorem{theorem}{Theorem}
\theoremstyle{remark}
\newtheorem{rem}{Remark}
\newtheorem{example}{Example}

\newcommand{\E}{\mathbf{E}}

\newcommand{\Rp}{\mathbb{R}_+}
\newcommand{\Zp}{\mathbb{Z}_+}
\newcommand{\I}{\textbf{I}}
\newcommand{\II}{\mathcal{I}}
\newcommand{\df}[1]{\mathcal{D}_{#1}}
\renewcommand{\vec}[1]{\text{\boldmath $#1$}}

\newcommand{\st}[3]{\bar{#1}_{#2}^{#3}}
\newcommand{\pst}[3]{\hat{#1}_{#2}^{#3}}
\newcommand{\T}[2]{T_{#1}^{#2}}

\parskip \smallskipamount

\begin{document}

 \begin{center}
   \Large\bfseries A note on insensitivity in stochastic networks
 \end{center}

 \begin{center}
   Stan Zachary
   \footnotetext{
     {\it American Mathematical Society 1991 subject classifications.\/}
     Primary 60K20


     }
 \end{center}

 \begin{center}
  \textit{Maxwell Institute for Mathematical Sciences, Heriot-Watt University\\
     Edinburgh}
 \end{center}

\begin{quotation}\small
  We give a simple and direct treatment of insensitivity in stochastic
  networks which is quite general and which provides probabilistic
  insight into the phenomenon.  In the case of multi-class networks,
  the results generalise those of Bonald and Prouti\`{e}re (2002,
  2003).
\end{quotation}

\section{Introduction}
  
It is well-known that many stochastic networks---notably queueing and
loss networks---have stationary distributions of their level of
occupancy which depend on certain input distributions only through the
means of the latter.  This phenomenon of \emph{insensitivity} has been
studied by various authors over an extended period of time, in varying
degrees of generality and abstraction, and using a variety of
techniques.

In the present paper we revisit this topic to develop an insight of
Pechinkin (1983, 1987) to give a very simple and direct treatment of
insensitivity.  In particular the approach avoids those based on
brute-force calculations, the consideration of phase-type
distributions (Schassberger, 1978, Whittle, 1985, Bonald and
Prouti\`{e}re, 2002, 2003), or the use of quite complex machinery for
handling generalised semi-Markov processes (Burman, 1981,
Schassberger, 1986)---although such processes are implicit in the
current approach.  It further avoids assumptions about, for example,
continuity of distributions, necessary for some of the above
approaches, and also explicitly identifies the entire stationary
distributions of the networks concerned, showing that, where
insensitivity obtains, these stationary distributions have a
particularly simple and natural form.  Pechinkin used his insight,
which involves what is in effect a coupling argument together with
induction, to give probabilistic proofs of the insensitivity of a
number of single-class loss systems with state-dependent arrival
rates---results originally proved analytically by Sevastyanov (1957).
He also indicated the wider applicability of the approach in the
single-class case.  In the present paper we give a substantial
reformulation of the underlying idea, under more general conditions
and showing that its most natural expression is in terms of balance
equations.  This considerably simplifies its application to
single-class systems---notably the quite complex coupling
constructions are no longer needed.  It further makes possible the
extension of the idea to the multi-class networks considered in
Section~\ref{sec:multi-class-networks}.  The main aim is to provide
probabilistic insight, notably for multi-class networks.  Indeed it is
shown that insensitivity is simply a byproduct, under appropriate
conditions, of probabilistic independence.

We study networks in which individuals arrive at various
\emph{classes} at rates which may depend on the state of the entire
system, bringing \emph{workloads} which are independent and
identically distributed within classes and which have finite means.
Within each class workloads are reduced at rates which may again be
state-dependent (when the rate is constant workloads may be identified
with lifetimes in classes), and on completion of its workload an
individual moves to a different class or leaves the system, with
probabilities which may yet again be state-dependent.

In order to obtain insensitivity we typically require that an
individual joining a class is immediately served, i.e.\ has its
workload reduced, at a rate which is the same as that of an individual
immediately prior to leaving the class (where in each of these cases
the number of individuals in each class of the system is the
same)---more generally that the service discipline should define a
network which is \emph{symmetric} in the sense of Kelly (1979).  The
most common example is that of processor-sharing networks, but other
possibilities are well-known, for example, ``last-in-first-out
preemptive resume'' networks.  We shall concentrate on a very broad
class of processor-sharing networks, introduced by Bonald and
Prouti\`{e}re (2002) and including, for example, traditional loss
networks and processor-sharing Whittle and Jackson networks, as
special cases).  We shall also indicate the simple modifications
required to deal with other possibilities.

For the above class of processor-sharing networks, Bonald and
Prouti\`{e}re used phase-type arguments to show that, under conditions
which correspond to the satisfaction of the appropriate partial
balance equations, the stationary distribution of the number of
individuals in each class is insensitive to the workload
distributions, subject to the means of the latter being fixed and to
the distributions themselves being drawn from the broad class of
Cox-type distributions (dense in the class of all distributions on
$\Rp$).  In the present paper we formally consider all workload
distributions on $\Rp$ with finite means, and identify also the
stationary residual workload distributions.  However, as stated above
our main aim is to give a direct and probabilistically natural
treatment.  It turns out (and is in many cases well-known) that, when
the appropriate partial balance equations are satisfied for such a
network, then the stationary distribution of the entire system,
\emph{including the specification of residual workloads}, is such that
departures from each class are exactly balanced by arrivals to that
class---in a sense again to be made precise below.  Indeed, for
single-class systems, this is the essence of Pechinkin's insight.
What is of interest is that same idea extends to establish
insensitivity for the very much more general networks considered here,
and indeed appears also to establish insensitivity in more abstract
settings such as that considered by Whittle (1985), though we do not
formally consider this more abstract environment here.

In order to fix ideas, it is convenient to consider first, in
Section~\ref{sec:single-class-networks}, single-class networks.  Here
the extension of previous ideas is not too difficult.  Nevertheless it
is desirable to give a careful treatment of this case, avoiding
notational complexity while preserving rigour, so as both to establish
the underlying principle and also to set the scene for the multi-class
networks which we consider in Section~\ref{sec:multi-class-networks}.

\section{Single-class networks}
\label{sec:single-class-networks}

Consider an open system with a single class of individual (customer,
call, or job).  Individuals arrive as a Poisson process with
state-dependent rate~$\alpha(n)$, where $n$ is the number of
individuals currently in the system.  Arriving individuals have
\emph{workloads} which are independent of each other and of the
arrivals process with a common distribution~$\mu$ on $\Rp$
which we assume to have a finite mean~$m(\mu)$.  
While there are $n$ individuals in the system, their \emph{total}
workload is reduced at a rate $\beta(n)\ge0$, where we assume
$\beta(n)>0$ if and only if $n>0$; an individual departs the system
when its workload is reduced to zero.  By suitably redefining the
rates~$\beta(n)$ if necessary, we may, and do, assume without loss of
generality that the mean workload $m(\mu)=1$.

We consider first the processor-sharing case.  Here when there are
$n>0$ individuals in the system, the workload of each is
simultaneously reduced at a rate $\beta(n)/n$, and the set-up
described above becomes a fairly general description of a single-class
processor sharing system.  A special case is the simple Erlang loss
system, in which, for some $\alpha,\beta>0$, we have
$\alpha(n)=\alpha\I(n<C)$ for some \emph{capacity}~$C\le\infty$ (where
$\I$ is the indicator function) and $\beta(n)=n\beta$ for all $n\ge0$.
Here individuals are typically referred to as calls, and workloads
correspond to call durations (since $\beta(n)/n$ is independent of
$n$).  A further special case is the $M/GI/m/\infty$ processor-sharing
queue, in which, again for some $\alpha,\beta>0$, we have
$\alpha(n)=\alpha$ for all $n$ and $\beta(n)=\min(n,m)\beta$ for all
$n$ and some fixed $m$.

We represent the system as a Markov process~$(X(t))_{t\ge0}$ by
defining its state at any time~$t$ to be the number~$n$ of individuals
then in the system together with their residual workloads at that
time.  (An alternative is to record, for each individual, the workload
completed at time~$t$.)  For given $n>0$ these workloads
form an (unordered) set, and may be regarded as taking
values in the quotient space~$S_n$ obtained from $\Rp^n$ by
identifying points which may be obtained from each other under
permutation of their coordinates.  The $\sigma$-algebra
$\mathcal{B}(S_n)$ on $S_n$ is similarly formed in the obvious manner
from the Borel $\sigma$-algebra on $\Rp^n$.  The state
space~$S$ for the process~$(X(t))_{t\ge0}$ is then the union of the
$S_n$, $n\ge0$, where the set $S_0$ is taken to consist of a single
point, and its associated $\sigma$-algebra~$\mathcal{B}(S)$ consists
of those sets which are countable unions of sets in the
$\sigma$-algebras~$\mathcal{B}(S_n)$.  The process~$(X(t))_{t\ge0}$ is
thus an instance of a piecewise-deterministic Markov process (Davis,
1984, 1993).  However, we avoid the need for most of the general
machinery for handling such processes.

We define the probability distribution~$\st{\mu}{}{}$ on $\Rp$ to be
the stationary residual life distribution of the renewal process with
inter-event distribution~$\mu$, that is, if $\mu$ has distribution
function~$F$ then $\st{\mu}{}{}$ has distribution function~$G$ given
by
\begin{displaymath}
  G(x) = 1 - \int_x^\infty (1-F(y))\,dy
\end{displaymath}
(recall that $m(\mu)=1$).  Note that the ``residual life'' here should
be thought of as a residual workload rather than a time.  For each
$n\ge1$, define also the probability distribution~$\st{\mu}{n}{}$ on
$S_n$ to be the product of $n$ copies of the
distribution~$\st{\mu}{}{}$, again with the above identification of
points in $\Rp^n$ (more formally,
$\st{\mu}{n}{}(A)=\st{\mu}{}{n}(\theta^{-1}(A))$, $A \in
\mathcal{B}(S_n)$, where $\st{\mu}{}{n}$ is the product of $n$ copies
of the distribution~$\st{\mu}{}{}$ and $\theta$ is the projection from
$\Rp^n$ into the quotient space~$S_n$).  Thus $\st{\mu}{n}{}$
represents the joint distribution of the residual lives at any time in
a set of $n$ independent stationary renewal processes each with
inter-event distribution~$\mu$; we define also $\st{\mu}{0}{}$ to be
the probability distribution concentrated on the single-point
set~$S_0$.  For each $n$, we also regard $\st{\mu}{n}{}$ as a
distribution on $S$, assigning its total mass one to the set~$S_n$.
Finally, for any distribution~$\pi$ on $\Zp$, define the distribution
$\st{\mu}{\pi}{}$ on $S$ by
$\st{\mu}{\pi}{}=\sum_{n\in\Zp}\pi(n)\st{\mu}{n}{}$.  Thus
$\st{\mu}{\pi}{}$ assigns probability~$\pi(n)$ to the event that there
are $n$ individuals in the system, and, conditional on this event,
assigns the distribution~$\st{\mu}{n}{}$ to their residual
workloads.

\begin{theorem}\label{thm:single}
  Suppose that the distribution~$\pi$ on $\Zp$ is the solution of the
  balance equations
  \begin{equation}
    \label{eq:1}
    \pi(n+1)\beta(n+1) = \pi(n)\alpha(n),
    \qquad n\ge0,
  \end{equation}
  and that 
  \begin{equation}
    \label{eq:2}
    \sum_{n\ge0}\pi(n)\alpha(n)<\infty.
  \end{equation}
  Then the distribution $\st{\mu}{\pi}{}$ on $S$ is stationary for the
  process~$(X(t))_{t\ge0}$, and in particular the distribution~$\pi$
  is stationary for the associated number of individuals in the
  system.
\end{theorem}

\begin{rem}
  The condition~\eqref{eq:2} ensures that, under stationarity,
  individuals arrive at the system at a finite rate.

\end{rem}

\begin{proof}[Proof of Theorem~\ref{thm:single}]
  In order to exclude pathological behaviour in the argument below, we
  make the one additional assumption that the distribution~$\mu$ has
  no atom of probability at zero.  This is without loss of generality:
  in the case that $\mu$ does have such an atom, the evolution of the
  system may clearly be equivalently described by redefining $\alpha$,
  $\beta$ and $\mu$ so as to remove it, and the result of the theorem
  is easily obtained via this reparametrisation.

  Analogously to the definition of $\st{\mu}{n}{}$, for each $n\ge1$,
  define the probability distribution~$\pst{\mu}{n}{}$ on $S_n$ to be
  the product of $n-1$ copies of the distribution~$\st{\mu}{}{}$ and a
  single copy of the distribution~$\mu$, yet again with the above
  identification of points in $\Rp^n$.  (More formally,
  $\pst{\mu}{n}{}(A)=\pst{\mu}{}{(n)}(\theta^{-1}(A))$,
  $A\in\mathcal{B}(S_n)$, where
  $\pst{\mu}{}{(n)}=\st{\mu}{}{n-1}\times\mu$ and $\theta$ is again
  the projection from $\Rp^n$ into the quotient space~$S_n$.)  We
  again regard $\pst{\mu}{n}{}$ as a distribution on $S$, assigning
  mass one to the set~$S_n$.

  Consider now the modified process~$(\hat{X}(t))_{t\ge0}$ on $S$
  describing the system in which the workload distribution is again
  $\mu$ and in which, when there are $n\ge1$ individuals in the
  system, individual workloads are again reduced at rate $\beta(n)/n$;
  however, for the modified system, (a)~an individual departing on
  completion of its workload is immediately replaced by another
  bringing an independent workload with distribution~$\mu$,
  (b)~external arrivals to the system are not accepted.  Thus, for the
  modified system, the number of individuals remains constant, and
  conditional on this being $n$, the system behaves as a set of $n$
  independent renewal processes, each of which has stationary residual
  workload distribution~$\st{\mu}{}{}$.  Hence, for any distribution
  $\pi'$ on $\Zp$, the distribution~$\st{\mu}{\pi'}{}$ on $S$ is
  stationary for the process~$(\hat{X}(t))_{t\ge0}$.

  Let $(P_t)_{t\ge0}$ and $(\hat{P}_t)_{t\ge0}$ be the semigroups of
  transition kernels associated respectively with the
  processes~$(X(t))_{t\ge0}$ and $(\hat{X}(t))_{t\ge0}$.  For any
  $a>0$, let $\df{a}$ be the class of functions~$f$ on $S$ taking
  values in $[0,1]$ and satisfying the continuity condition
  \begin{equation}
    \label{eq:3}
    \left\lvert (P_t f(x) - f(x)) \right\rvert \le at
    \qquad\text{for all $x\in S$ and $t>0$},
  \end{equation}
  where $P_tf(x)=\int_S P_t(x,dy)f(y)$.  For any such $f$ and for any
  distribution~$\nu$ on $S$, define also $\nu f = \int_S f(x) \nu(dx)$
  and, for any $t>0$, define $\nu{}P_tf=\nu(P_t{}f)$ (so that
  $\nu{P_t}f$ is the expectation of $f(X(t))$ when $(X(t))_{t\ge0}$ is
  given initial distribution~$\nu$); similarly define $\nu\hat{P_t}f$.

  Now compare the behaviour of the processes~$(X(t))_{t\ge0}$ and
  $(\hat{X}(t))_{t\ge0}$, each started with the distribution
  $\st{\mu}{\pi}{}$; so as to simplify the description below we couple
  these two processes so that they agree until the time of the first
  arrival or workload completion.  We then have (see the further
  explanation below) that, with this common initial distribution, for
  any $a>0$, $f\in\df{a}$, and $h>0$,
  \begin{align}
    \st{\mu}{\pi}{} P_h f - \st{\mu}{\pi}{} \hat{P}_h f
    & = \E \bigl(f(X(h)) - f(\hat{X}(h))\bigr)\nonumber \\
    & = h \sum_{n\ge0} \pi(n)
    \left[
      \alpha(n)(\pst{\mu}{n+1}{}f - \st{\mu}{n}{}f)
      + \beta(n)(\st{\mu}{n-1}{}f - \pst{\mu}{n}{}f)
    \right]
    + o(h)    \label{eq:4}\\
    & = h \sum_{n\ge0} 
    \left[
      \pi(n) \alpha(n) - \pi(n+1) \beta(n+1)
    \right]
    (\pst{\mu}{n+1}{}f - \st{\mu}{n}{}f)
    + o(h)    \nonumber\\
    & = o(h)  \label{eq:5}
  \end{align}
  as $h\to0$ (recall in \eqref{eq:4} that $\beta(0)=0$); further the
  above convergence as $h\to0$ is uniform over $f\in\df{a}$ in the
  sense that (\ref{eq:5}) may be written as
  \begin{equation}
    \label{eq:6}
    \sup_{f\in\df{a}}
    \lvert
    \st{\mu}{\pi}{} P_h f - \st{\mu}{\pi}{} \hat{P}_h f
    \rvert
    = o(h) \qquad\text{as $h\to0$}
  \end{equation}
  (again see below).  To show~\eqref{eq:4} note first that, from the
  above coupling and for any $h>0$, we have $f(X(h))=f(\hat{X}(h))$
  except where there is either at least one external arrival or at
  least one workload completion in $[0,h]$.  It follows from the
  definition of $\st{\mu}{\pi}{}$ that, conditional on
  the number of individuals initially being $n$ the probability of an
  external arrival in $[0,h]$ is $\alpha(n)h+o(h)$ as $h\to0$, and
  that an arriving individual finding the distribution of the system
  to be $\st{\mu}{n}{}$ changes this to $\pst{\mu}{n+1}{}$ in the case
  of the process~$(X(t))_{t\ge0}$ and leaves it unchanged in the case
  of the process~$(\hat{X}(t))_{t\ge0}$.  Similarly, again conditional
  on the number of individuals initially being $n$ (and recalling that
  $m(\mu)=1$), the probability of a workload completion in a time
  interval~$[0,h]$ is $\beta(n)h+o(h)$ as $h\to0$, and that under the
  distribution~$\st{\mu}{n}{}$, conditional on such a completion
  taking place, the residual workload distribution becomes
  $\st{\mu}{n-1}{}$ in the case of the process~$(X(t))_{t\ge0}$ and
  $\pst{\mu}{n}{}$ in the in the case of the
  process~$(\hat{X}(t))_{t\ge0}$.  Further it follows from the
  conditions~\eqref{eq:1} and \eqref{eq:2} that, under the initial
  distribution~$\st{\mu}{\pi}{}$, the probability of two or more
  arrivals or workload completions in $[0,h]$ is $o(h)$ as $h\to0$.
  That the relation~\eqref{eq:4} now holds as $h\to0$ with the
  uniformity over $f\in\df{a}$ required for \eqref{eq:6} follows
  easily from these results and from the definition of $\df{a}$.  To
  see this note that, since $f\in\df{a}$ implies that $f$ takes values
  in $[0,1]$, the contribution to the error term in (\ref{eq:4})
  resulting from the neglect of the possibility of two or more
  arrivals or workload completions in $[0,h]$ is uniformly $o(h)$ as
  $h\to0$ as required.  Similarly the terms
  $\pst{\mu}{n+1}{}f-\st{\mu}{n}{}f$ and
  $\st{\mu}{n-1}{}f-\pst{\mu}{n}{}f$ in (\ref{eq:4}) are obtained by
  treating the precise time of the first arrival or workload
  completion within $[0,h]$ as if it were time $h$;  (recalling
  that $\pst{\mu}{n+1}{}$, etc, are probability measures) it follows
  from \eqref{eq:3} that the consequent error in each of the above
  two terms is bounded by $2ah$, so that the further contribution
  to the error term in (\ref{eq:4}) is $O(h^2)$ as $h\to0$, again with
  uniformity over $f\in\df{a}$. The
  relations~\eqref{eq:5}, and hence (\ref{eq:6}), are now immediate
  from the balance equations~\eqref{eq:1}.  Since the
  distribution~$\st{\mu}{\pi}{}$ is stationary for the
  process~$(\hat{X}(t))_{t\ge0}$, it now follows from~\eqref{eq:6}
  that, again for any $a>0$ and $h>0$,
  \begin{displaymath}
    \sup_{f\in\df{a}}
    \lvert
    \st{\mu}{\pi}{} P_h f - \st{\mu}{\pi}{} f
    \rvert
    = o(h)
   \qquad\text{as $h\to0$}.
  \end{displaymath}
  Further, it is straightforward that if $f\in\df{a}$, then also
  $P_tf\in\df{a}$ for any $t>0$.  Standard manipulations using the
  semigroup structure of $(P_t)_{t\ge0}$, e.g.\ the consideration of
  increasingly refined partitions of the interval~$[0,t]$, now give
  that, for all $a>0$, $f\in\df{a}$, and $t\ge0$,
  \begin{equation}\label{eq:7}
    \st{\mu}{\pi}{} P_t f = \st{\mu}{\pi}{} f.
  \end{equation}
  Finally, we show that it follows from (\ref{eq:7}) that
  $\st{\mu}{\pi}{}P_t=\st{\mu}{\pi}{}$ for all $t>0$, so that
  $\st{\mu}{\pi}{}$ is stationary for $(X(t))_{t\ge0}$ as required.
  It is sufficient to show that, for any $n\ge1$ and any
  set~$A\in\mathcal{B}(S_n)$ whose inverse image in
  $\Rp^n$ under the mapping $\theta$ defined above is a product of
  intervals in $\Rp$, we have
  \begin{equation}
    \label{eq:8}
    \st{\mu}{\pi}{} P_t \I_A = \st{\mu}{\pi}{} \I_A,
  \end{equation}
  where $\I_A$ is the indicator function of the set~$A$.  It follows
  from the piecewise deterministic form of the
  process~$(X(t))_{t\ge0}$ that we may choose a sequence of
  functions~$(f_k,\,k\ge1)$ such that, for each $k$, (i)
  $f_k\in\df{a}$ for some $a>0$ and (ii) $f_k$ and $\I_A$ agree except
  on a set whose Lebesgue measure (under $\theta^{-1}$) in $\Rp^n$
  tends to zero as $k\to\infty$.  Since $\st{\mu}{\pi}{}$, and so also
  $\st{\mu}{\pi}{}P_t$, are non-atomic distributions, the
  result~(\ref{eq:8}) now follows by using (\ref{eq:7}) with $f=f_k$
  and letting $k\to\infty$.
\end{proof}

\begin{rem}\label{rem:balance}
  Suppose that the equations~\eqref{eq:1} above are multiplied by the
  signed measure $(\pst{\mu}{n+1}{}-\st{\mu}{n}{})$ to give
  \begin{equation}
    \label{eq:9}
    \pi(n)\alpha(n)(\pst{\mu}{n+1}{} - \st{\mu}{n}{})
    = \pi(n+1)\beta(n+1)(\pst{\mu}{n+1}{} - \st{\mu}{n}{}),
    \qquad n\ge0.
  \end{equation}
  These equations have an obvious interpretation as representing,
  under the distribution~$\st{\mu}{\pi}{}$ and for each $n\ge0$, a detailed
  balance of flux between $S_n$ and $S_{n+1}$, not just with regard to
  the total probability assigned to each of these spaces, but also with
  regard to the distribution of the residual workload sizes: the
  intuition underlying the derivation of \eqref{eq:5} above---which is
  also that of Pechinkin's coupling approach---is that, under
  $\st{\mu}{\pi}{}$, an arrival finding $n$ individuals in the system
  transforms the residual workload distribution from $\st{\mu}{n}{}$ to
  $\pst{\mu}{n+1}{}$, while a departure from the system when it
  contains $n+1$ individuals transforms the residual workload
  distribution from what would have been $\pst{\mu}{n+1}{}$, if the
  individual had remained in the system with a renewed workload, to
  the distribution $\st{\mu}{n}{}$.
\end{rem}

In the case where we do not have processor-sharing, i.e.\ in which it
is no longer the case that at any time all workloads are being reduced
at the same rate, it is necessary at any time to distinguish the
individuals in the system.  Thus each~$S_n$ above is replaced by
$\Rp^n$ and the state space $S$ is replaced by
$S^*=\bigcup_{n\ge0}\Rp^n$.  We consider as an example the case of
the single-server queue with ``last-in-first-out preemptive resume''
discipline, in which at any time all service effort is devoted to the
last individual to arrive at the system.  If at any time there are $n$
individuals in the system, we may index these by $i=1,\dots,n$ in the
order of their arrival, and no individual changes index during its
time in the system; as usual arrivals occur as a Poisson process with
rate $\alpha(n)$, and the workload of individual~$n$ is now being
reduced at rate~$\beta(n)$, while that of the remaining individuals is
being reduced at rate~$0$.  As previously, define the probability
distribution~$\st{\mu}{}{}$ on $\Rp$ to be the stationary residual
life distribution of the renewal process with inter-event
distribution~$\mu$, and, for each $n\ge0$, let the distribution
$\st{\mu}{n}{}$ on $\Rp^n$ be now the (ordered) product of $n$ copies
of $\st{\mu}{}{}$.  For each $n\ge1$, let the
distribution~$\pst{\mu}{n}{}$ on $\Rp^n$ be the (ordered) product of
$n-1$ copies of the distribution~$\st{\mu}{}{}$ and a single copy of
the distribution~$\mu$, with the latter assigned to the $n$th
coordinate of $\Rp^n$.  Finally, for any distribution~$\pi$ on $\Zp$,
define the distribution $\st{\mu}{\pi}{}$ on $S$ by
$\st{\mu}{\pi}{}=\sum_{n\in\Zp}\pi(n)\st{\mu}{n}{}$ as before.  With
these (re)definitions, both Theorem~\ref{thm:single} and its proof
remain unchanged as stated.  Again the underlying reason is as given
in Remark~\ref{rem:balance} above: under the
distribution~$\st{\mu}{\pi}{}$, and relative to the modified process
considered in the proof of Theorem~\ref{thm:single}, an arrival
finding $n$ individuals in the system transforms the residual workload
distribution from $\st{\mu}{n}{}$ to $\pst{\mu}{n+1}{}$, while a
departure from the system when it contains $n+1$ individuals
transforms the residual workload distribution from $\pst{\mu}{n+1}{}$
to $\st{\mu}{n}{}$.  Note that this balance does not obtain in the
case of, for example, a ``first-in-first-out'' discipline, and here,
as is again well known, we do not have the above insensitivity.

\section{Multi-class networks}
\label{sec:multi-class-networks}

Consider now a multi-class network.  We concentrate on the
processor-sharing case---adaptations to other disciplines may be made
as in the single-class case.  Let $\II=\{1,\dots,N\}$ denote the set of
classes, and let $\vec{n}=(n_i,\,i\in\II)$ where $n_i$ is the number
of individuals in each class~$i$.  An individual entering class~$i$
acquires a workload which has distribution~$\mu^i$ with nonzero finite
mean~$m(\mu^i)$; we again assume without loss of generality that
$m(\mu^i)=1$; the workload of each individual in class $i$ is reduced
at a state-dependent rate~$\phi_i(\vec{n})/n_i$, where
$\phi_i(\vec{n})>0$ if and only if $n_i>0$.  Individuals arrive at
each class~$i$ from outside the network as a Poisson process with
state-dependent rate~$\phi_{0i}(\vec{n})$; on completion of its
workload in any class~$i$ an individual moves to class~$j$ with
state-dependent probability~$\phi_{ij}(\vec{n})/\phi_i(\vec{n})$ or
leaves the network with
probability~$\phi_{i0}(\vec{n})/\phi_i(\vec{n})$, where
\begin{equation}\label{eq:10}
  \sum_{j\in\II}\phi_{ij}(\vec{n}) + \phi_{i0}(\vec{n}) = \phi_i(\vec{n})
\end{equation}
(there are no problems in allowing the possibility
$\phi_{ii}(\vec{n})>0$).  The workloads, arrivals processes and routing
decisions are all independent.

As in the single-class case, we represent the system as a Markov
process~$(X(t))_{t\ge0}$ by defining its state at any time to be the
vector~$\vec{n}$ introduced above together with the residual workloads
at that time of the set of individuals in each class.  For given
$\vec{n}$, these workloads take values in the space~$S_\vec{n}$ which
is the ordered product of the spaces~$S_{n_i}$, $i\in\II$, where, as
previously, each $S_{n_i}$ is formed from $\Rp^{n_i}$ by identifying
points within the latter space which may be obtained from each other
under permutation of their coordinates (and where again the set~$S_0$
contains a single point).  The state space~$S$ for the system is the
union of all the possible $S_\vec{n}$, and the spaces $S_\vec{n}$ and
$S$ are endowed with the obvious
$\sigma$-algebras~$\mathcal{B}(S_\vec{n})$ and $\mathcal{B}(S)$.  

Analogously to the single-class case, for each $i\in\II$, define
the probability distribution~$\st{\mu}{}{i}$ on $\Rp$ to be the
stationary residual life distribution of the renewal process with
inter-event distribution~$\mu^i$ (as previously the residual life
should be interpreted as a residual workload).  For each $i\in\II$ and
for each $n_i\ge1$, define as previously the
distribution~$\st{\mu}{n_i}{i}$ on $S_{n_i}$ to be the product of
$n_i$ copies of the distribution~$\st{\mu}{}{i}$ (again with the above
identification of points in $\Rp^{n_i}$)---representing the joint
distribution of the residual lives in a set of $n_i$ independent
stationary renewal processes each with inter-event
distribution~$\mu^i$; define also $\st{\mu}{0}{i}$ to be the
probability distribution concentrated on the single-point set~$S_0$.
For each $\vec{n}\in\Zp^N$, define the distribution
$\st{\vec{\mu}}{\vec{n}}{}$ on $S_\vec{n}$ to be the (ordered) product
distribution which, for each $i\in\II$, assigns the
distribution~$\st{\mu}{n_i}{i}$ to $S_{n_i}$.  We again regard the
distribution~$\st{\vec{\mu}}{\vec{n}}{}$ as a distribution on $S$,
assigning its total mass one to the set~$S_\vec{n}$.  For any positive
distribution~$\pi$ on $\Zp^N$, we define the distribution
$\st{\vec{\mu}}{\pi}{}$ on $S$ by
$\st{\vec{\mu}}{\pi}{}=\sum_{\vec{n}\in\Zp^N}\pi(\vec{n})\st{\vec{\mu}}{\vec{n}}{}$.

It is notationally convenient to expand the set $\II$ to
$\II'=\{0\}\cup\II$, treating $0$ as an extra class feeding external
arrivals to, and receiving departures from, the network.  (However,
the components of the state~$\vec{n}$ of the network remain indexed in
the original set~$\II$.)  For completeness we define
$\phi_{00}(\vec{n})=0$ for all $\vec{n}$, and also
$\phi_{ij}(\vec{n})=0$ for all $i\in\II$, $j\in\II'$, and $\vec{n}$
such that $n_i=0$ (so that \eqref{eq:10} above remains valid for such
$\vec{n}$ also).  For each $i\in\II$, let $\vec{e}_i$ be the
$N$-dimensional vector whose $i$th component is $1$ and whose other
components are $0$, and let $\vec{e}_0$ be the $N$-dimensional vector
all of whose components are $0$.  For each $\vec{n}$ and each
$i,j\in\II'$ define the vector
$\T{i}{j}\vec{n}=\vec{n}-\vec{e}_i+\vec{e}_j$; define also
$T_i\vec{n}=\T{i}{0}\vec{n}=\vec{n}-\vec{e}_i$ and
$T^j\vec{n}=\T{0}{j}\vec{n}=\vec{n}+\vec{e}_j$.

\begin{theorem}\label{thm:multi}
  Suppose that the distribution~$\pi$ on $\Zp^N$ satisfies the partial
  balance equations
  \begin{equation}
    \label{eq:11}
      \pi(\vec{n})\sum_{j\in\II'}\phi_{ij}(\vec{n})
      = \sum_{j\in\II'}\pi(\T{i}{j}\vec{n})\phi_{ji}(\T{i}{j}\vec{n}),
      \qquad
      \vec{n}\in\Zp^N, \quad i\in\II',
  \end{equation}
  where, for $\vec{n}$ and $i\in\II$ such that $n_i=0$ we interpret
  the right side of \eqref{eq:11} as zero (recall that when $n_i=0$ we
  have $\phi_{ij}(\vec{n})=0$ for all $j\in\II'$ so that \eqref{eq:11}
  is automatically satisfied in this case).  Suppose also that
  \begin{equation}
    \label{eq:12}
    \sum_{\vec{n}\in\Zp^N}\pi(\vec{n})\sum_{i\in\II}\phi_{0i}(\vec{n})<\infty.
  \end{equation}
  Then $\st{\vec{\mu}}{\pi}{}$ is stationary for the
  process~$(X(t))_{t\ge0}$, and in particular $\pi$ is stationary for
  the associated numbers of individuals in the system.  Conversely, if
  a distribution~$\pi$ on $\Zp^N$ is stationary for the numbers of
  individuals in the system for all $\vec{\mu}=(\mu^i,\,i\in\II)$ such
  that $m(\mu^i)=1$ for all $i\in\II$, then $\pi$ satisfies the
  equations~(\ref{eq:11}).
\end{theorem}

\begin{proof}
  Suppose first that $\pi$ satisfies the equations~\eqref{eq:11}.  As
  in the proof of Theorem~~\ref{thm:single}, we again assume without
  loss of generality that each distribution~$\mu^i$ has no atom of
  probability at zero.

  For each $\vec{n}$ and for each $i$ such that $n_i\ge1$, define also
  the residual workload distribution~$\pst{\vec{\mu}}{\vec{n}}{i}$ on
  $S_\vec{n}$ by
  $\pst{\vec{\mu}}{\vec{n}}{i}%
  =\pst{\mu}{n_i}{i}\times\prod_{j\ne{}i}\st{\mu}{n_j}{j}$
  where $\pst{\mu}{n_i}{i}$ is defined as in the proof of
  Theorem~\ref{thm:single}.  Thus $\pst{\vec{\mu}}{\vec{n}}{i}$
  corresponds to each individual in each class~$j$ having
  independently the stationary residual workload
  distribution~$\st{\mu}{}{j}$, except only that a single individual
  in the class $i$ is given the workload distribution~$\mu^i$.  For
  each $\vec{n}$, define also
  $\pst{\vec{\mu}}{\vec{n}}{0}=\st{\vec{\mu}}{\vec{n}}{}$.

  Again as in the proof of Theorem~\ref{thm:single}, define the
  process~$(\hat{X}(t))_{t\ge0}$ on $S$ to be that appropriate to the
  modified system in which there are no arrivals, departures, or
  transfers between classes; rather each individual in each class~$i$,
  on completion of its workload, acquires a new independent workload
  with distribution~$\mu^i$.  Thus the occupancy of the system remains
  constant; conditional on this being $\vec{n}$, individual workloads
  in any class~$i$ such that $n_i>1$ are again reduced at rate
  $\phi_i(\vec{n})/n_i$ and the system behaves as a set of independent
  renewal processes.  Further, for any distribution $\pi'$ on $\Zp^N$,
  the distribution~$\st{\vec{\mu}}{\pi'}{}$ on $S$ is stationary for
  $(\hat{X}(t))_{t\ge0}$.

  Again let $(P_t)_{t\ge0}$ and $(\hat{P}_t)_{t\ge0}$ be the
  semigroups of transition kernels associated respectively with the
  processes~$(X(t))_{t\ge0}$ and $(\hat{X}(t))_{t\ge0}$, and, for any
  $a>0$, let $\df{a}$ be the class of functions~$f$ on $S$ taking
  values in $[0,1]$ and satisfying the earlier continuity
  condition~\eqref{eq:3}.  Comparison of the behaviour of the
  processes~$(X(t))_{t\ge0}$ and $(\hat{X}(t))_{t\ge0}$, each started
  with the distribution $\st{\vec{\mu}}{\pi}{}$ and coupled as in the
  earlier proof until the time of the first external arrival or
  workload completion, now gives that, for any $a>0$, $f\in\df{a}$,
  and $h>0$,
  \begin{align}
    \st{\vec{\mu}}{\pi}{} P_h f - \st{\vec{\mu}}{\pi}{} \hat{P}_h f
    \hspace{-6em} & \nonumber \\[1ex]
    & = \E \bigl(f(X(h)) - f(\hat{X}(h))\bigr)\nonumber \\
    & = h \sum_{\vec{n}\in\Zp^N}\pi(\vec{n})\sum_{i\in\II'}\sum_{j\in\II'}
    \phi_{ij}(\vec{n})
    \left(
      \pst{\vec{\mu}}{\T{i}{j}\vec{n}}{j}f - \pst{\vec{\mu}}{\vec{n}}{i}f
    \right)
    + o(h)    \label{eq:13}\\
    & = h \sum_{\vec{n}\in\Zp^N}
    \Biggl(
    \sum_{i\in\II'}\sum_{j\in\II'}\pi(\vec{n})
    \phi_{ji}(\vec{n})
    \pst{\vec{\mu}}{\T{j}{i}\vec{n}}{i}f
    - \sum_{i\in\II'}\sum_{j\in\II'}\pi(\vec{n})
    \phi_{ij}(\vec{n})
    \pst{\vec{\mu}}{\vec{n}}{i}f 
    \Biggr)
    + o(h)    \nonumber \\
    & = h \sum_{\vec{n}\in\Zp^N}
    \Biggl(
    \sum_{i\in\II'}\sum_{j\in\II'}\pi(\T{i}{j}\vec{n})
    \phi_{ji}(\T{i}{j}\vec{n})
    \pst{\vec{\mu}}{\vec{n}}{i}f
    - \sum_{i\in\II'}\sum_{j\in\II'}\pi(\vec{n})
    \phi_{ij}(\vec{n})
    \pst{\vec{\mu}}{\vec{n}}{i}f 
    \Biggr)
    + o(h)    \nonumber \\
    & = h \sum_{\vec{n}\in\Zp^N}\sum_{i\in\II'}
    \Biggl(
    \sum_{j\in\II'}\pi(\T{i}{j}\vec{n})
    \phi_{ji}(\T{i}{j}\vec{n})
    - \sum_{j\in\II'}\pi(\vec{n})
    \phi_{ij}(\vec{n})
    \Biggr)
    \pst{\vec{\mu}}{\vec{n}}{i}f
    + o(h)    \nonumber \\
    & = o(h) \label{eq:14}
  \end{align}
  as $h\to0$, with uniformity of convergence over all $f\in\df{a}$, so
  that we may write
  \begin{equation}
    \label{eq:15}
    \sup_{f\in\df{a}}
    \lvert
    \st{\vec{\mu}}{\pi}{} P_h f - \st{\vec{\mu}}{\pi}{} \hat{P}_h f
    \rvert
    = o(h) \qquad\text{as $h\to0$}
  \end{equation}
  (recall again that $\phi_{ij}(\vec{n})=0$ whenever $n_i=0$, so that
  there is no difficulty with the lack of a formal definition of
  $\pst{\vec{\mu}}{\vec{n}}{i}$ in this case).  The
  identity~(\ref{eq:15}) is
  simply the multi-class version of the identity~\eqref{eq:6} in the
  proof of Theorem~\ref{thm:single}, and is similarly obtained, albeit
  with a slightly more compact notation; in particular, conditional on
  the common initial distribution of the two processes being given by
  $\st{\vec{\mu}}{\vec{n}}{}$, in the time interval $[0,h]$ where $h$
  is small, a transition from $i$ to $j$ in the original
  system---where either $i$ or $j$ may be $0$---occurs with
  probability $\pi(\vec{n})\phi_{ij}(\vec{n})h+o(h)$ as $h\to0$, and
  in this case the distribution of the process~$(X(t))_{t\ge0}$
  becomes $\pst{\vec{\mu}}{\T{i}{j}\vec{n}}{j}$ while that of the
  process~$(\hat{X}(t))_{t\ge0}$ becomes
  $\pst{\vec{\mu}}{\vec{n}}{i}$; that \eqref{eq:13} now holds with the
  required uniformity of convergence follows, using also \eqref{eq:12},
  as in the earlier proof; finally the results \eqref{eq:14}, and so
  also (\ref{eq:15}), follow
  from the partial balance equations~\eqref{eq:11}.  Since the
  distribution~$\st{\vec{\mu}}{\pi}{}$ is stationary for the
  process~$(\hat{X}(t))_{t\ge0}$, it now follows from~\eqref{eq:15}
  that, again for any $a>0$ and $h>0$,
  \begin{displaymath}
    \sup_{f\in\df{a}}
    \lvert
    \st{\vec{\mu}}{\pi}{} P_h f - \st{\vec{\mu}}{\pi}{} f
    \rvert
    = o(h)
    \qquad\text{as $h\to0$}.
  \end{displaymath}
  That $\st{\vec{\mu}}{\pi}{}$ is now stationary for $(X(t))_{t\ge0}$
  follows as in the proof of Theorem~\ref{thm:single}.

  Now suppose that a distribution~$\pi$ on $\Zp^N$ is stationary for
  the numbers of individuals in the system for all
  $\vec{\mu}=(\mu^i,\,i\in\II)$ with $m(\mu^i)=1$ for all $i\in\II$.
  A proof that $\pi$ then necessarily satisfies the partial balance
  equations~\eqref{eq:11} is given by Bonald and Prouti\`{e}re (2002,
  2003).  In summary, consider the case in which in every class the
  workload distribution is exponential with mean~$1$, and, for any
  fixed class~$i$, compare this with the case in which, for some
  $0<\lambda<1$, the workload distribution in class~$i$ is replaced by
  a mixture of two distributions, obtained by choosing with
  probability~$\lambda$ an exponential distribution with
  mean~$\lambda^{-1}$, and with probability~$1-\lambda$ the
  distribution concentrated on $0$.  Both these models may be
  (re)formulated as simple Markov jump processes---in the latter case
  the transition rates into and out of the class~$i$ are reduced by a
  factor~$\lambda$.  Since $\pi$ is stationary in both cases,
  comparison of the (full) balance equations for stationarity yields
  the partial balance equations~\eqref{eq:11}.
\end{proof}

\begin{example}
  \emph{Processor-sharing Whittle networks.}  Suppose that for some
  $\nu\ge0$, some strictly positive function~$\Phi$ on $\Zp^N$, and
  some stochastic matrix~$P=(p_{ij},\,i\in\II',\,j\in\II')$ such that
  $p_{00}=0$, we have, for each $\vec{n}$,
  \begin{alignat*}{2}
    \phi_{0j}(\vec{n}) & = \nu p_{0j}, & \qquad & j \in \II', \\
    \phi_{ij}(\vec{n}) & = \frac{\Phi(T_i\vec{n})}{\Phi(\vec{n})}\, p_{ij},
    && i \in \II, \quad j \in \II',
  \end{alignat*}
  where we again make the convention that $\Phi(T_i\vec{n})=0$
  whenever $n_i=0$.  Then it is readily checked that the partial
  balance equations~\eqref{eq:11} are satisfied by
  \begin{equation}\label{eq:16}
    \pi(\vec{n}) = a \Phi(\vec{n}) \prod_{i\in\II} \rho_i^{n_i},
  \end{equation}
  for any $a>0$ and positive solution $\vec{\rho}=(\rho_i,\,i\in\II)$ of the
  equations
  \begin{align*}
    \nu & = \sum_{j\in\II}\rho_j p_{j0},\\
    \rho_i & = \sum_{j\in\II}\rho_j p_{ji} + \nu p_{0i}, \qquad i\in\II.
  \end{align*}
  Thus in particular the stationary distribution~$\pi$ given by
  \eqref{eq:16} for the number of individuals of each type in the
  system is insensitive to the $\mu^i$ (recall our assumption
  $m(\mu^i)=1$ for all $i$).  For the case where $P$ is irreducible
  and $\nu>0$, the above equations for $\vec{\rho}$ have a unique
  solution.  Again when $P$ is irreducible and when $\nu=0$
  (corresponding to a closed network) $\pi$ remains uniquely
  determined, up to a multiplicative constant, by \eqref{eq:16}.  The
  case where $\Phi(\vec{n})=\prod_{i\in\II}\lambda_i^{n_i}$ for
  positive constants $(\lambda_i,\,i\in\II)$ characterises
  processor-sharing Jackson networks.  Further discussion of Whittle
  networks is given by Serfozo (1999) and, for processor-sharing
  networks, by Bonald and Prouti\`{e}re (2002, 2003).
\end{example}

\begin{example}
  \emph{Networks with no internal transitions.}  Suppose that
  $\phi_{ij}(\vec{n})=0$ for all $i,j\in\II$ and for all
  $\vec{n}\in\Zp^N$, so that no transitions are possible between the
  classes in $\II$.  The partial balance equations~\eqref{eq:11} then
  reduce to the detailed balance equations
  \begin{equation}
    \label{eq:17}
    \pi(\vec{n})\phi_{i0}(\vec{n})
    = \pi(T_i\vec{n})\phi_{0i}(T_i\vec{n}),
      \qquad
      \vec{n}\in\Zp^N, \quad n_i\ge1, \quad i\in\II.
  \end{equation}
  (In the case of a single class, these equations further reduce to
  the equations~\eqref{eq:1}.)  An example is given by a traditional
  (uncontrolled) loss network---see, for example, Kelly (1986).  This
  is naturally processor-sharing.  Here workloads are identified with
  call durations and, for some set $\mathcal{A}\subset\Zp^N$ such that
  $\vec{n}\in\mathcal{A}$ implies $T_i\vec{n}\in\mathcal{A}$ for all
  $\vec{n}$ and $i$ such that $n_i\ge1$ ($\mathcal{A}$ is typically
  defined by capacity constraints), we have
  \begin{align*}
    \phi_{0i}(\vec{n}) & = \nu_i\I(T^i\vec{n}\in\mathcal{A})\\
    \phi_{i0}(\vec{n}) & = \sigma_i n_i,
  \end{align*}
  for some vectors $(\nu_i,\,i\in\II)$ and $(\sigma_i,\,i\in\II)$ of
  strictly positive parameters.  The equations~\eqref{eq:17} are then
  satisfied by
  \begin{equation}
    \label{eq:18}
        \pi(\vec{n}) = a \prod_{i\in\II} \frac{\kappa_i^{n_i}}{n_i!},
  \end{equation}
  where $\kappa_i=\nu_i/\sigma_i$ for each $i$, and where $a$ is
  naturally chosen to be a normalising constant.  As was originally
  shown by Burman \textit{et al} (1984), we therefore again have
  insensitivity of the occupancy distribution~$\pi$ of the network.
  The stationary distribution of the residual call durations is as
  identified by Theorem~\ref{thm:multi}.

  Other examples of processor-sharing networks with no internal
  transitions are given by those used to model connections in
  communications networks with simultaneous resource requirements and
  variable bandwidth requirements---see, for example, Bonald and
  Massouli\'e (2001) and de Veciana \textit{et al} (2001).  Here it is
  far from automatic that the detailed balance equations~\eqref{eq:17}
  are satisfied.
\end{example}

\section*{Acknowledgements}

The author is most grateful to Serguei Foss, Takis Konstantopoulos and
Ilze Ziedins for some very helpful discussions, and to the referee for
a careful reading of the manuscript and some most helpful suggestions.

\end{document}